\documentclass{amsart}

\usepackage{graphicx} 
\usepackage{amsmath}
\usepackage{amsthm}
\numberwithin{equation}{section} 
\usepackage{amssymb}
\usepackage{empheq}
\usepackage{cases}
\usepackage{verbatim}
\newtheorem{theorem}{Theorem}[section]
\newtheorem{proposition}{Proposition}[section]
\newtheorem{lemma}{Lemma}[section]
\newtheorem{corollary}{Corollary}[section]
\newtheorem{definition}{Definition}
\usepackage{caption}
\usepackage{subcaption}
\theoremstyle{remark}
\newtheorem{remark}{Remark}[section]
\usepackage[margin=1in]{geometry}
\usepackage[T1]{fontenc}    
\usepackage{lmodern}  
\makeatletter
\renewcommand\subsubsection{\@startsection{subsubsection}{3}%
  \z@{.5\linespacing\@plus.7\linespacing}{-.5em}%
  {\normalfont\bfseries}}
\makeatother

\title{Inverted parabola as a stable steady solution of the Surface Growth Model}
\author{Won-Suk Lee}
\address{Department of Physics \& Astronomy, Seoul National University, 1 Gwanak-ro, Gwanak-gu, Seoul, 08826, Republic of Korea}
\email{lookatme@snu.ac.kr}


\begin{document}
\begin{abstract}
    For the surface growth model under molecular beam epitaxy on $\mathbb{R}$, we show that an inverted parabola is a stable steady solution in the sense that any small perturbation of it decays in $L^\infty(\mathbb{R})$ at the rate $O(t^{-1/2}\log t)$. The proof is based on the explicit solution operator generated by the linearized perturbation equation. This result partially justifies the merging of two neighboring islands into a bigger one. Moreover, any derivative of its solution is shown to be decaying in $L^2(\mathbb{R})$ for sufficiently small initial data.
\end{abstract}

\maketitle


\section{Introduction}

In this paper we consider the initial value problem of molecular beam epitaxy (MBE) in surface growth model (SGM). Though it is motivated by the 2-dimensional physics, we consider the 1-dimensional problem
\begin{equation}\label{eq: SGM}
\begin{split}
    \left\{ \begin{array}{ccc}
    u_t + u_{xxxx}+ u_{xx} &= -(u_x^2)_{xx}, \\
    u(0,x)&=u_0(x),
    \end{array}\right.
\end{split}
\end{equation}
where $x\in\mathbb{R}$, $t\ge0$. This equation is to be solved for the time-dependent surface height $u(t,x)$ from the reference level. It is based on the continuity equation where the time derivative term is in general proportional to the Laplacian of the chemical potential terms and that explains why every term with spatial derivatives includes second derivatives. The second derivative term $u_{xx}$ reflecting the linear instability is sometimes omitted, but in this paper we include this term into the analysis. The nonlinear term $-(u_x^2)_{xx}$ is often called conserved KPZ term which is introduced by Villain \cite{villain_1991} and contributes to the breaking of up-down symmetry of the system. It is related to the phenomena observed in the numerical simulation \cite{STEIN_WINKLER_2005} (and also in experimental results \cite{Thin_Film_Materials_Freund}).

In recent decades, several mathematicians have studied stability using a Lyapunov functional (e.g. \cite{Blomker_Romito_2009}) or in the sense of statistics (e.g. \cite{Blomker_Flandoli_Romito_2009}), as well as well-posedness and regularity (e.g. \cite{Blomker_Romito_2009}, \cite{Blomker_Romito_2012}, \cite{Cheng_Li_Wei_2026}, \cite{Ozanski_Robinson_2019} and \cite{STEIN_WINKLER_2005}) . Also, it has been known that there is a quadratic steady solution $u(t,x)=a-(1/4)(x-b)^2$(e.g. \cite{2000_Raible_Linz_Hanggi}). Numerical results point out that the solution smoothly grows except for the points of local minima which are considered to be the candidates for singularities (for details, see \cite{STEIN_WINKLER_2005}). However, finite-time singularity formation from general initial data is still an open problem for \eqref{eq: SGM}.

There are three types of steady solutions and among them the inverted parabolae often appear in numerical simulations mentioned above. While the details are described in the appendix (summarized in Theorem \ref{thm: classification of steady solutions}), it is supported by the physics in the background of the surface growth model. The stationarity is obtained by balancing diffusion and uphill currents that stop particles from spreading out. Chemical potential makes particles prefer places with lower curvature. Consequently it becomes unfavorable to form a steep hill. Meanwhile, there also exist uphill currents --- e.g. one due to Schwoebel barrier --- that prevent the collapse of hills. It can be compared to a sand grain failing to slip down from the slope due to friction. In such natural situations where there are competing opposite forces and the system has a preferred direction, a concave curve often turns out to be an equilibrium and even stable.

As the initial data $u_0(x)$ evolves to $u(t,x)$, we define the perturbation $v(t,x)=u(t,x)-[(1/2)-(1/4)x^2]$, whose behavior is governed by \eqref{eq: perturbed SGM}. Our goal is to prove the asymptotic stability of $v(t,x)=0$ (corresponding to the steady solution itself) and to determine the rate at which it converges to zero. Here $v_0$ represents an initial perturbation from the quadratic steady part, so we have to impose restrictions on the scale of $v_0$. That condition is stated in Theorem \ref{thm: main result}.
\begin{equation}\label{eq: perturbed SGM}
    v_t+v_{xxxx}-(xv_{xx})_x=-\partial_{x}^2(v_x^2).
\end{equation}

\begin{theorem}\label{thm: main result}
    Let the initial data $v_0$ be in $(L^1\cap H^5)(\mathbb{R})$. Then there exists $\epsilon>0$ such that the condition $||v_0||_{H^5(\mathbb{R})}\le\epsilon$ is sufficient for the existence and uniqueness of a classical global-in-time solution $v(t,x)$ of \eqref{eq: perturbed SGM} converging to zero in $L^\infty(\mathbb{R})$. If $v_0\in H^\sigma(\mathbb{R})$ with an integer $\sigma>5$, then this regularity persists. Also, $||v(t,\cdot)||_{L^\infty(\mathbb{R})}=O(t^{-1/2}\log t)$.
\end{theorem}

\begin{remark}
    $v_0$ being small is decided by the following conditions:
    \begin{itemize}
        \item $\gamma(\epsilon_0)<1$, where $\gamma$ is defined in Theorem \ref{thm: bootstrap closure for small initial data} of Section \ref{sec: stability of infinite parabola}.
        \item For some real number $\rho>1/2$, it is satisfied that 
        \begin{equation}
            ||v_0||_{L^2(\mathbb{R})}\le\min_{1\le n\le5}\left\{  \frac{1}{C_n}\sqrt{\frac{2\rho-1}{\rho^2(\rho+1)}(2n-1)} \right\},
        \end{equation}
        where $C_n$ is a sequence of constant numbers for Moser-type estimate in \eqref{eq: CS and Moser for V_3}.
        \item $||\partial_x^mv_0||_{L^2(\mathbb{R})}\le\sqrt{2\rho^m}\epsilon_0$ for $1\le m \le 5$ and $||v_0||_{L^2(\mathbb{R})}\le\epsilon_0$.
    \end{itemize}
\end{remark}

\begin{remark}
    Requirement of $L^1$ is in accordance with the finite mass deviation from the steady solution $(1/2)-(1/4)x^2$, though it does not exactly match due to the issue of sign. This interpretation can be extended to the 2-dimensional problem, from which the problem is motivated.
\end{remark}

\begin{remark}
    In short, small deviation from the parabola disappears. If a curve with two hills whose maximal points are sufficiently close --- e.g. a small narrow Gaussian perturbation --- then they will merge into a single bigger hill. Figure \ref{fig: L^infty time evolution of v(t,x)} visualizes the decay of $v(t,x)$ with the Gaussian initial data. Though it is not about the detailed mechanism, this supports the plausibility of the coarsening process generally known in material science. For more details, see \cite{Thin_Film_Materials_Freund}.
\end{remark}

\begin{figure}
    \centering
    \begin{subfigure}{0.53\textwidth}
        \includegraphics[width=\linewidth]{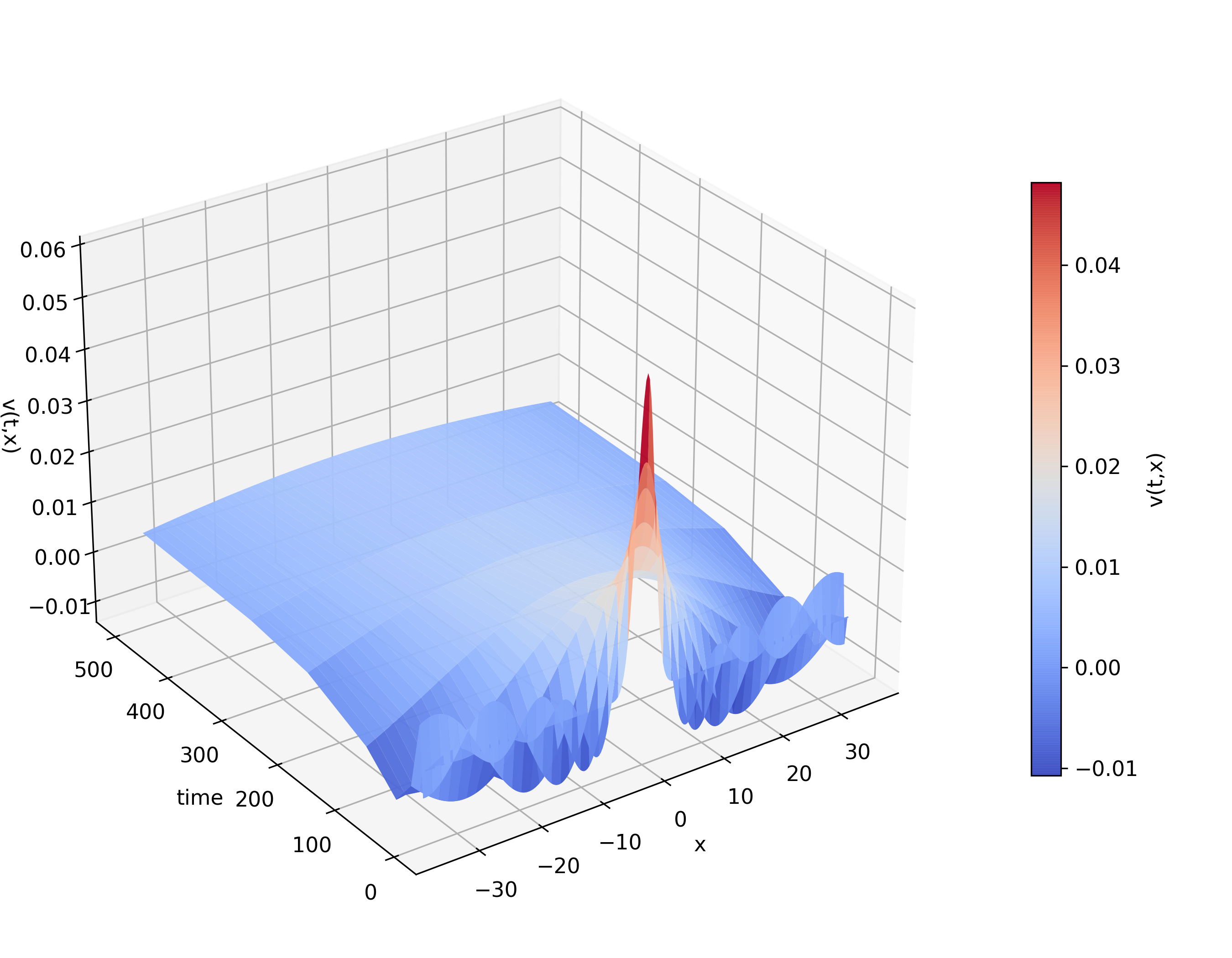}
        \caption{Time evolution of $v(t,x)$. Even though the initial data is given as a small Gaussian function, the solution almost immediately becomes oscillatory for some time.}
        \label{fig: time evolution of perturbation profile}
    \end{subfigure}
    \hfill
    \begin{subfigure}{0.45\textwidth}
        \includegraphics[width=\linewidth]{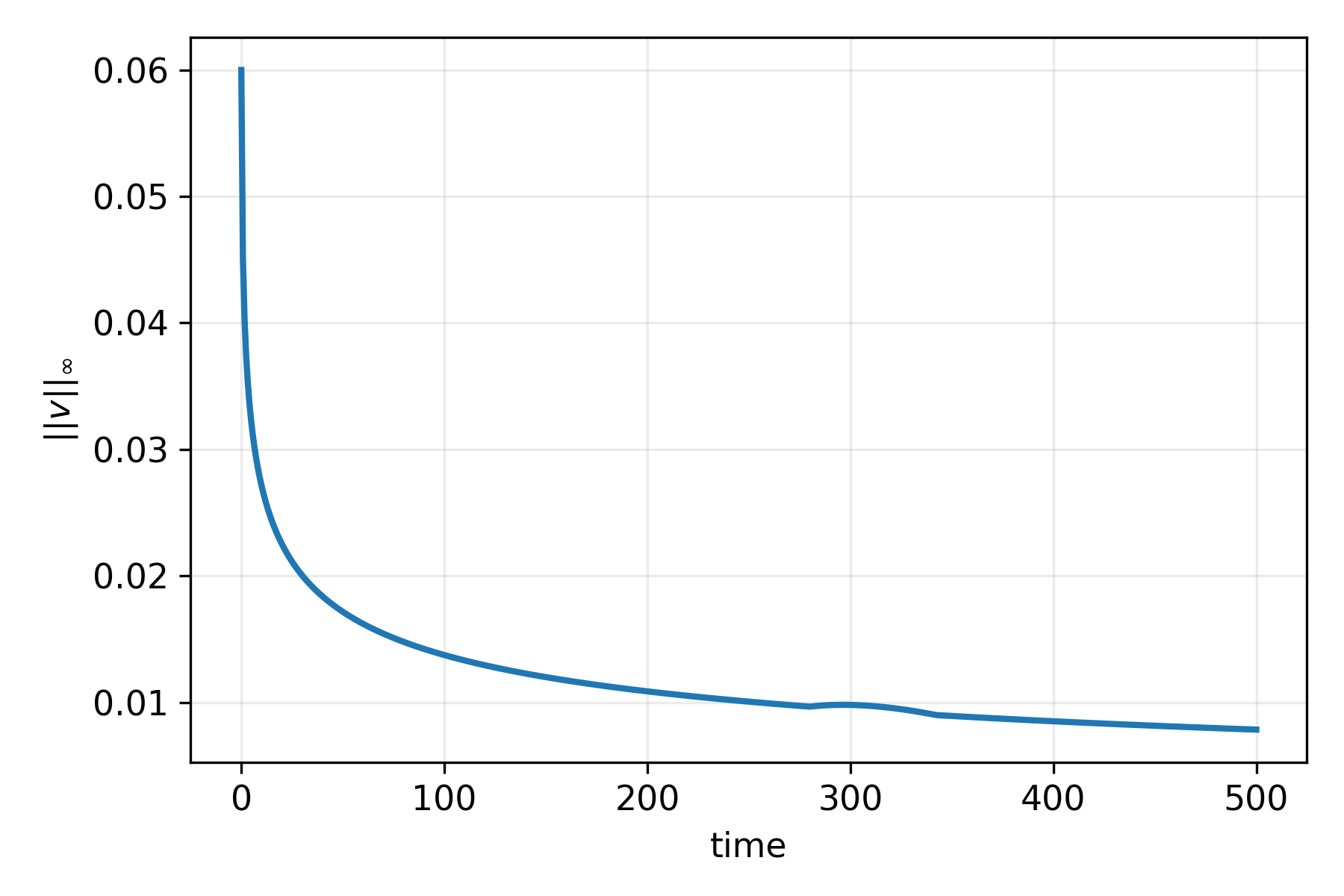}
        \caption{The decay of $||v(t,\cdot)||_{L^\infty(\mathbb{R})}$. One can see that it decreases to zero, but it sometimes increases due to the spatial oscillation of $v(t,x)$ shown in (a).}
        \label{fig: decay of L^infty of v(t,x)}
    \end{subfigure}

    \caption{Time evolution of the perturbation $v(t,x)$ and its $L^\infty$ norm.}
    \label{fig: L^infty time evolution of v(t,x)}
\end{figure}

The remaining sections are organized as follows: Before considering \eqref{eq: perturbed SGM}, we ignore the nonlinear term and analyze the linearized equation in Section \ref{sec. linear perturbation}. Some stability properties and additional regularities given by the solution-generating semigroup related to \eqref{eq: linearized perturbed SGM} will also be discussed. In Section \ref{sec: L2 estimates of the perturbation} we get back to the equation including the nonlinear term and try to obtain some a priori energy estimates in $L^2(\mathbb{R})$. Based on those results, we prove the stability of \eqref{eq: perturbed SGM} in the sense of Theorem \ref{thm: main result} in Section \ref{sec: stability of infinite parabola}. The well-posedness of \eqref{eq: perturbed SGM} is discussed in Section \ref{sec: well-posedness}, to support the validity of all the preceding statements.


\section{Linearized Perturbation of SGM}\label{sec. linear perturbation}
Before considering the full equation for $v(t,x)$, we linearize the equation around $v=0$,

\begin{equation}\label{eq: linearized perturbed SGM}
    w_t + w_{xxxx}-(xw_{xx})_{x} = 0.
\end{equation}

As this solution can be explicitly solved, we can construct many useful statements directly from the formula. In addition to the following proposition, in Section \ref{sec: S(t) regularities} we investigate $L^2$ estimates related to \eqref{eq: linearized perturbed SGM}.

\begin{proposition}\label{prop: GWP of lin pert}
    If the initial data $w_0\in(L^1\cap H^1)(\mathbb{R})$ is given, then there is a solution $w(t,x)$ to \eqref{eq: linearized perturbed SGM} whose explicit form is expressed in \eqref{eq: explicit linear perturbation}. Moreover, for any nonnegative integer $n$, $||\partial_x^nw(t,\cdot)||_{L^\infty(\mathbb{R})}$ satisfies the decay $O(t^{-(n+1)/2}\log t)$ as $t\rightarrow\infty$.
\end{proposition}

Throughout this paper the Fourier transform works as one of the most important tools for calculations. We adopt the following unitary convention of the Fourier transform and Sobolev norms.

\begin{definition}
    For a function $f:\mathbb{R}\rightarrow\mathbb{R}$, Fourier transform $\mathcal{F}[f]$ is defined as below and equivalently denoted by $\hat{f}$.
    \begin{equation}
       \mathcal{F}[f](\xi):=\frac{1}{\sqrt{2\pi}}\int_\mathbb{R}f(x)e^{-i\xi x}dx.
    \end{equation}
    As a consequence, Fourier inverse transform is given as follows.
    \begin{equation}
        \mathcal{F}^{-1}[\hat{f}](x):=\frac{1}{\sqrt{2\pi}}\int_\mathbb{R}\hat{f}(\xi)e^{i\xi x}d\xi.
    \end{equation}
\end{definition}

\begin{definition}
    A function $f:\mathbb{R}\rightarrow\mathbb{R}$ that can be Fourier transformed is in Sobolev space $H^\sigma$ if the following quantity is finite. Also, we use this quantity as the norm of $H^\sigma(\mathbb{R})$.
    \begin{equation}
        ||f||_{H^\sigma(\mathbb{R})}:=\left(\int_\mathbb{R}(1+\xi^2)^\sigma|\hat{f}(\xi)|^2d\xi\right)^{\frac{1}{2}}.
    \end{equation}
\end{definition}

\subsection{Explicit Solution of the Linearized Equation}\label{subsec: explicit formula for linear perturbation}
After executing the Fourier transform $\hat{w}:=\mathcal{F}[w]$, \eqref{eq: linearized perturbed SGM} becomes the transport equation in the frequency domain.

\begin{equation}\label{eq: FT of lin pert SGM}
    \hat{w}_t+(k^4+k^2)\hat{w}-\partial_k(k^3\hat{w})=0.
\end{equation}

Since the first derivative in $k$ indicates that \eqref{eq: FT of lin pert SGM} is transport type, we look for the trajectory map $\phi(t)$ such that $f(t):=\hat{w}(t,\phi(t))$ satisfies:

\begin{equation}
    \frac{df}{dt}=\hat{w}_t-k^3 \hat{w}_k.
\end{equation}
This leads to the ODE for the trajectory: $\phi'(t)=-\phi(t)^3$, which can be explicitly solved with the given initial data $\phi(0)=k_0$.

\begin{equation}\label{eq: characteristic curve for lin pert}
    \phi(t)=\frac{k_0}{\sqrt{1+2k_0^2t}}.
\end{equation}
It is notable that $\phi(t)$ is bounded in $(-1/\sqrt{2t},1/\sqrt{2t})$, even though the initial data is distributed all over the frequency space. That is, $\hat{w}(t,\cdot)$ is compactly supported for all $t>0$ and thus $w(t,\cdot)$ is analytic according to the Paley-Wiener theorem. Along the characteristic curve $(t,\phi(t))$ lying in the $t-k$ plane, the solution for \eqref{eq: FT of lin pert SGM} is

\begin{equation}
    \hat{w}(t,\phi(t))=f(t) = f(0)\exp\left(-\int_0^t(\phi(s)^4-2\phi(s)^2)ds\right). \\
\end{equation}
Here $f(0)=\hat{w}(0,k(0))=\hat{w}_0(k_0)$ and $\phi(s)$ is given by \eqref{eq: characteristic curve for lin pert}. Also, the integrations in the exponents are calculated as follows.

\begin{equation}
\begin{split}
    \int_0^t\phi(s)^4ds &= \int_0^t\frac{k_0^4}{(1+2k_0^2s)^2}ds = \frac{k_0^4t}{1+2k_0^2t}, \\
    \int_0^t2\phi(s)^2ds &= \int_0^t\frac{2k_0^2}{1+2k_0^2s}ds = \log(1+2k_0^2t). \\
\end{split}
\end{equation}

Now we can head to the solution of \eqref{eq: linearized perturbed SGM} with the inverse Fourier transform of $\hat{w}(t,k)$. Substituting $k$ with $\phi(t)$ and fixing $t$, we can express the inverse Fourier transform as an integration over $k_0$. It is easy to check that it actually fulfills \eqref{eq: linearized perturbed SGM} and the initial condition.

\begin{equation}\label{eq: explicit linear perturbation}
\begin{split}
    w(t,x) &= \frac{1}{\sqrt{2\pi}}\int_{-\infty}^\infty \hat{w}(t,k)e^{ikx}dk \\
    &= \frac{1}{\sqrt{2\pi}}\int_{-\infty}^\infty\hat{w}_0(k_0)(1+2k_0^2t)^{-1/2}\exp\left(-\frac{k_0^4t}{1+2k_0^2t}+\frac{ik_0x}{\sqrt{1+2k_0^2t}}\right)dk_0.
\end{split}
\end{equation}

\subsection{Stability of Linearized Perturbation}\label{sec: regularity of lin pert}
Now we estimate the $L^\infty$ norm of the linear perturbation $w$. This shows the asymptotic stability of \eqref{eq: linearized perturbed SGM} near the steady parabolic solution and in the later sections it is a cornerstone for proving the stability of nonlinear perturbations. We start from the absolute value of the perturbation $w(t,x)$. It is also assumed that $t>e$.

\begin{equation}
    |w(t,x)|\le\frac{1}{\sqrt{2\pi}}\int_{-\infty}^\infty|\hat{w}_0(k_0)|(1+2k_0^2t)^{-1/2}\exp\left(-\frac{k_0^4t}{1+2k_0^2t}\right)dk_0.
\end{equation}
Separating the domain of integration into I: $k_0^2t<1$ and II: $k_0^2t\ge1$,

\begin{equation}
\begin{split}
    &\int_I (1+2k_0^2t)^{-1/2}\exp\left(-\frac{k_0^4t}{1+2k_0^2t}\right)dk_0 \le \int_Idk_0=\frac{2}{\sqrt{t}}, \\
    &\int_{II} (1+2k_0^2t)^{-1/2}\exp\left(-\frac{k_0^4t}{1+2k_0^2t}\right)dk_0\le \int_{t^{-1/2}}^{\infty}\sqrt{\frac{2}{t}}\frac{1}{k_0}e^{-k_0^2/3}dk_0 \\
    &\le \sqrt{\frac{2}{t}}\left(\int_{t^{-1/2}}^{1}\frac{1}{k_0}dk_0 + \frac{1}{2}\int_{-\infty}^{\infty}e^{-k_0^2/3}dk_0\right)\le \sqrt{\frac{1}{2t}}\left(\log t +\sqrt{3\pi}\right). \\
\end{split}
\end{equation}
As a result, 

\begin{equation}
    |w(t,x)|\lesssim ||\hat{w}_0||_{L^\infty(\mathbb{R})}\frac{\log t}{\sqrt{t}}\le ||w_0||_{L^1}\frac{\log t}{\sqrt{t}}.
\end{equation}
where the last line comes from the fact that $\log t>1$. Note that this bound is especially valid for large $t$ and becomes singular near zero. However, since \eqref{eq: explicit linear perturbation} recovers the initial data when $t=0$, the continuity-in-time of $w(t,x)$ leads to the temporal boundedness of $w(t,x)$ in addition to the above asymptotic bound. In detail, we have the boundedness by $||\hat{w}_0||_{L^1(\mathbb{R})}$.

\begin{equation}
    |w(t,x)|\le\frac{1}{\sqrt{2\pi}}\int_{-\infty}^\infty|\hat{w}_0(k_0)|dk_0.
\end{equation}

$||\hat{w}_0||_{L^1(\mathbb{R})}$ can be also bounded by the quantities from the real space. We temporarily set the free parameter $R$ that splits the frequency domain and will later use it to optimize the estimate of $||\hat{w}_0||_{L^1(\mathbb{R})}$.

\begin{equation}
\begin{split}
    \int_{-\infty}^\infty|\hat{w}_0(k_0)|dk_0&\le \left(\int_{|k_0|<R}dk_0\right)^\frac{1}{2}||w_0||_{L^2(\mathbb{R})}+\left(\int_{|k_0|\ge R}\frac{dk_0}{k_0^2}\right)^\frac{1}{2}||\partial_x w_0||_{L^2(\mathbb{R})}\\
    &= (2R)^\frac{1}{2}||w_0||_{L^2(\mathbb{R})}+\left(\frac{2}{R}\right)^{\frac{1}{2}}||\partial_x w_0||_{L^2(\mathbb{R})}\le 2\sqrt{2}||w_0||_{H^1(\mathbb{R})},
\end{split}
\end{equation}
where we choose $R=||\partial_x w_0||_{L^2(\mathbb{R})}/||w_0||_{L^2(\mathbb{R})}$ for optimization and the Cauchy-Schwarz inequality has been adopted for the last inequality. Then the complete form of estimate for $||w(t,\cdot)||_{L^\infty(\mathbb{R})}$ valid for all $t\ge0$ is

\begin{equation}\label{eq: L-infty estimate of lin pert}
    ||w(t,\cdot)||_{L^\infty(\mathbb{R})} \le \alpha_0\left(||w_0||_{L^1(\mathbb{R})}+||w_0||_{H^1(\mathbb{R})}\right)\frac{\log (e+t)}{\sqrt{1+t}}.
\end{equation}

When the initial perturbation $w_0(x)$ is given as the Gaussian function $e^{-x^2}$, we look for the value of its solution $w(t,x)$ at $x=0$, the point of symmetry.

\begin{equation}
    w(t,0) = \frac{1}{\sqrt{4\pi}}\int_{-\infty}^{\infty}(1+2k_0^2t)^{-1/2}\exp\left(-k_0^2\left(\frac{1}{4}+\frac{k_0^2t}{1+2k_0^2t}\right)\right)dk_0.
\end{equation}
Here we observe that the exponent in the integrand is greater than $-3k_0^2/4$ and in the region II defined above $(1+2k_0^2t)^{-1/2}>(3k_0^2t)^{-1/2}$ is satisfied. Thus we have the lower bound of $||w(t,\cdot)||_{L^\infty(\mathbb{R})}$ in the same form with the upper bound.

\begin{equation}
    ||w(t,\cdot)||_{L^\infty(\mathbb{R})} \ge |w(t,0)|\ge \frac{1}{\sqrt{4\pi}} \int_{II}(3k_0^2t)^{-1/2}\exp\left(-\frac{3}{4}k_0^2\right)dk_0\ge \frac{c'\log t}{\sqrt{t}}.
\end{equation}
This implies that the upper bound \eqref{eq: L-infty estimate of lin pert} cannot be improved in general. Exactly the same argument for $\partial_x^nw(t,x)$ results in the upper bounds of $||\partial_x^nw(t,\cdot)||_{L^\infty(\mathbb{R})}$. For any non-negative integer $n$,

\begin{equation}
\begin{split}
    |\partial_x^nw(t,\cdot)| &\le \frac{1}{\sqrt{2\pi}}\int_{-\infty}^\infty|\hat{w}_0(k_0)|\frac{|k_0|^n} {(1+2k_0^2t)^{(1+n)/2}}\exp\left(-\frac{k_0^4t}{1+2k_0^2t}\right)dk_0 \\
    &\le \frac{||\hat{w}_0||_{L^\infty(\mathbb{R})}}{\sqrt{2\pi}}\left(\int_{I}|k_0|^ndk_0+(2t)^{-\frac{n+1}{2}}\int_{II}\frac{1}{|k_0|}e^{-\frac{k_0^2}{3}}dk_0\right) \\
    &\lesssim ||w_0||_{L^1(\mathbb{R})}t^{-\frac{n+1}{2}}\log t. \\
\end{split}
\end{equation}

In summary,

\begin{equation}\label{eq: L-infty estimate of n-th deriv lin pert}
    ||\partial_x^nw(t,\cdot)||_{L^\infty(\mathbb{R})} \le \alpha_n||w_0||_{L^1(\mathbb{R})}t^{-(n+1)/2}\log t.
\end{equation}

\subsection{Properties of $S(t)$}\label{sec: S(t) regularities}

The operator $S(t)$ which generates the solution of \eqref{eq: linearized perturbed SGM} satisfies inequalities useful for the later sections. In this section we list them and verify ones not proven above. While looking into the stability of $w(t,x)$ we have seen $L^\infty$ properties.

\begin{itemize}
    \item If $u\in(L^1\cap H^1)(\mathbb{R})$, then $S(t)u$ is bounded in time and converges to zero, uniformly. In particular, for the asymptotic stability of $w=0$, $u\in L^1$ is sufficient.
    \item If $u\in L^1(\mathbb{R})$, then $\partial_x^nS(t)u$ always converges to zero as $t\rightarrow\infty$ and the higher order derivatives decay faster.
\end{itemize}

But there are also some properties related to $L^2$. 

\begin{lemma}\label{lemma: L2 of S(t)u}
    If $u\in L^2(\mathbb{R})$, then $||S(t)u||_{L^2(\mathbb{R})}\le ||u||_{L^2(\mathbb{R})}(1+t)^{1/4}$.
\end{lemma}

\begin{lemma}\label{lemma: L2 of D^n[S(t)u]}
    If $u\in H^n(\mathbb{R})$ and $n>0$, then $||\partial_x^nS(t)u||_{L^2(\mathbb{R})}\le ||\partial_x^nu||_{L^2(\mathbb{R})}$.
\end{lemma}

\begin{lemma}\label{lemma: regularity of S(t)[D^alpha u]}
    Given that $u\in H^j(\mathbb{R})$ and $t<1$, then there is a positive constant $\mu_{\alpha \sigma j}$ determined by a nonnegative integer $\alpha$ and positive real numbers $\sigma\ge j$ such that:
    \begin{equation}\label{eq: regularity of S(t)[D^alpha u]}
        ||S(t)[\partial_x^\alpha u]||_{H^\sigma(\mathbb{R})}\le \mu_{\alpha \sigma j}t^{-\frac{\sigma-j+\alpha}{4}}||u||_{H^j(\mathbb{R})}.
    \end{equation}
\end{lemma}

\subsubsection{$||S(t)u||_{L^2(\mathbb{R})}$ may grow at the rate $t^{1/4}$}

Even though $S(t)u$ has been shown to be stable in $L^\infty$, it is not in $L^2$. We first apply Plancherel's theorem to \eqref{eq: explicit linear perturbation}.
\begin{equation}
\begin{split}
    ||S(t)u||_{L^2(\mathbb{R})}^2=\int_{\mathbb{R}} |\hat{u}(k_0)|^2(1+2k_0^2t)^{1/2}\exp\left(-\frac{2k_0^4t}{1+2k_0^2t}\right)dk_0.
\end{split}
\end{equation}

As we substitute $z:=2k_0^2t$, the integrand other than $\hat{u}$ becomes $(1+z)^{1/2}\exp(-z^2/2t(1+z))$, which is maximized at $z=(-2+t+\sqrt{4+t^2})/2$ or equivalently $t=z(z+2)/(z+1)$. At this point, the maximal value is
\begin{equation}
    (1+z)^{\frac{1}{2}}\exp\left(-\frac{z^2}{2z(z+2)}\right)\le \left(1+z+\frac{z}{z+1}\right)^{\frac{1}{2}}\cdot 1 = (1+t)^{\frac{1}{2}},
\end{equation}
and hence less than $(1+t)^{1/2}$. 
\begin{equation}\label{eq: time scale of lin pert in L2}
    ||S(t)u||_{L^2(\mathbb{R})}^2\le \int_{\mathbb{R}}|\hat{u}(k_0)|^2(1+t)^{\frac{1}{2}}dk_0\le ||u||_{L^2(\mathbb{R})}^2(1+t)^{\frac{1}{2}}.
\end{equation}

\subsubsection{$\partial_x^nS(t)u$ is bounded in $L^2$}

Unlike the previous case, any spatial derivative of $S(t)u$ is guaranteed to be bounded in $L^2$ if the same order derivative of the initial data $u$ is bounded in $L^2$. As it can be shown by the direct differentiation of \eqref{eq: explicit linear perturbation}, that is because the factor $(1+2k_0^2t)^{1/2-n}$ becomes less than 1 unlike the case $n=0$.

\begin{equation}
\begin{split}
    ||\partial_x^nS(t)u||_{L^2(\mathbb{R})}^2 &= \int_{\mathbb{R}}|\hat{u}(k_0)|^2 \frac{|k_0|^{2n}}{(1+2k_0^2t)^{n-1/2}}\exp\left(-\frac{2k_0^4t}{1+2k_0^2t}\right)dk_0 \\
    &\le \int_\mathbb{R} |k_0^n\hat{u}(k_0)|^{2}dk_0 = ||\partial_x^n u||_{L^2(\mathbb{R})}^2.
\end{split}
\end{equation}

\subsubsection{Additional regularity is obtained by $S(t)$}

For the last one, we separate the integration contained in $S(t)$ into three parts; $I_1:|k_0|<t^{-1/4}$, $I_2: t^{-1/4}\le |k_0|<t^{-1/2}$ and $I_3:|k_0|\ge t^{-1/2}$. For any $t\ge0$,

\begin{equation}
    ||S(t)[\partial_x^\alpha u]||_{H^\sigma(\mathbb{R})}^2= \int_\mathbb{R}\left(1+\frac{k_0^2}{1+2k_0^2t}\right)^{\sigma}|k_0^\alpha\hat{u}(k_0)|^2(1+2k_0^2t)^{\frac{1}{2}}\exp\left(-\frac{2k_0^4t}{1+2k_0^2t}\right)dk_0.
\end{equation}

For simplicity any factors other than $\hat{u}$ have been combined to be denoted as $M_t(k_0)$. Then our goal is to show that $M_t(k_0)/(1+k_0^2)^j$ is bounded by $t^{-(\sigma-j+\alpha)/2}$. When $|k_0|<t^{-1/4}$, we use that $k_0^2t<t^{1/2}<1$ because $t<1$. 

\begin{equation}
\begin{split}
    M_t(k_0)&:=\left(1+\frac{k_0^2}{1+2k_0^2t}\right)^\sigma|k_0|^{2\alpha}(1+2k_0^2t)^{\frac{1}{2}}\exp\left(-\frac{2k_0^4t}{1+2k_0^2t}\right) \\
    &\le \sqrt{3}(1+k_0^2)^\sigma |k_0|^{2\alpha} \le \sqrt{3}(1+k_0^2)^j\left(2t^{-\frac{1}{2}}\right)^{\sigma-j}t^{-\frac{\alpha}{2}}. \\
\end{split}
\end{equation}
Writing $\mu_{\alpha \sigma j}^{(1)}=\sqrt{3}\cdot 2^{\sigma-j}$ for simplicity,

\begin{equation}
    \frac{M_t(k_0)}{(1+k_0^2)^j} \le \mu_{\alpha\sigma j}^{(1)}t^{-\frac{\sigma-j+\alpha}{2}}.
\end{equation}

When $t^{-1/4}\le |k_0|<t^{-1/2}$, $k_0^2t<1$ is still valid and we introduce $-2k_0^4t/(1+2k_0^2t)\le -(2/3)k_0^4t$ for exponent.

\begin{equation}
    \frac{M_t(k_0)}{(1+k_0^2)^j}\le \mu_{\alpha\sigma j}^{(1)}t^{-\frac{\sigma-j+\alpha}{2}}(k_0^4t)^{\frac{\alpha}{2}}\exp\left(-\frac{2}{3}k_0^4t\right) \le \mu_{\alpha\sigma j}^{(2)}t^{-\frac{\sigma-j+\alpha}{2}}.
\end{equation}

For the last case $|k_0|\ge t^{-1/2}$, the exponential factor is less than $\exp(-k_0^2/2)$.

\begin{equation}\label{eq: additional regularity from S(t): the third part of domain}
    \frac{M_t(k_0)}{(1+k_0^2)^j}\le (1+k_0^2)^{\sigma-j}|k_0|^{2\alpha}(4k_0^2t)^{1/2}\exp\left(-\frac{k_0^2}{2}\right)\le \mu_{\alpha \sigma j}^{(3)}t^{\frac{1}{2}}.
\end{equation}

As \eqref{eq: additional regularity from S(t): the third part of domain} offers better decay than $t^{-(\sigma-j+\alpha)/2}$, \eqref{eq: regularity of S(t)[D^alpha u]} is derived from summing up all those results. The positive constant $\mu_{\alpha\sigma j}$ is defined as $\mu_{\alpha\sigma j}^2:=\max_{i}\mu_{\alpha\sigma j}^{(i)}$.

\begin{equation}
    ||S(t)[\partial_x^\alpha u]||_{H^\sigma(\mathbb{R})}^2 \le \int_{\mathbb{R}}\mu_{\alpha sj}^2t^{-\frac{\sigma-j+\alpha}{2}}(1+k_0^2)^j|\hat{u}(k_0)|^2dk_0 =\mu_{\alpha\sigma j}^2t^{-\frac{\sigma-j+\alpha}{2}}||u||_{H^j(\mathbb{R})}^2.
\end{equation}


\section{Energy Estimate Including Nonlinearity}\label{sec: L2 estimates of the perturbation}

As we have looked for some $L^2$ properties of $S(t)$ in the previous section, we try to obtain the similar kind of estimate of the nonlinear perturbation $v(t,x)$ in $L^2$. We define $E_n(t):=(1/2)\int_{\mathbb{R}}|\partial_{x}^{n} v(t,x)|^2dx$ for nonnegative integers $n$. In the rest of this section we omit the domain and the variable of integral since they are the whole real line and the spatial variable $x$ in every line. 

\begin{lemma}
    For each $n\ge1$, if the perturbation \eqref{eq: perturbed SGM} satisfies $||v_{x}(t,\cdot)||_{L^\infty(\mathbb{R})}<\sqrt{(2-1/\rho)(2n-1)}/C_n$, then $E_{n}(t)$ is nonincreasing, as in \eqref{eq: L2 evol of E_n by higher terms}.
\end{lemma}

By substituting \eqref{eq: perturbed SGM} into the time derivative of $E_n$, it can be expressed in terms of higher order terms $E_{n+1}$ and $E_{n+2}$.

\begin{equation}
    \frac{d}{dt}E_n = \int(\partial_x^nv)\partial_x^n(\partial_tv)= \int (\partial_x^nv)(-\partial_x^{n+4}v+\partial_x^{n+1}(xv_{xx})-\partial_x^{n+2}(v_x^2)).
\end{equation}

For those three terms, we denote them as $V_1$, $V_2$, and $V_3$, respectively.

\begin{equation}
    V_1 = \int-(\partial_x^nv)(\partial_x^{n+4}v)=-\int|\partial_x^{n+2}v|^2 = -2E_{n+2}.
\end{equation}

\begin{equation}
\begin{split}
    V_2 &= \int(\partial_x^nv)\partial_x^{n+1}(xv_{xx})= -\int(\partial_x^{n+1}v)(n\partial_x^{n+1}v+x\partial_x^{n+2}v) \\
    &=-2nE_{n+1}-\frac{1}{2}\int x\partial_x\left(|\partial_x^{n+1}v|^2\right)=-(2n-1)E_{n+1}.
\end{split}
\end{equation}

\begin{equation}\label{eq: CS and Moser for V_3}
    V_3 = \int(\partial_x^nv)(\partial_x^{n+2}(v_x^2))\le ||\partial_x^{n+2}v||_{L^2(\mathbb{R})}||\partial_x^n(v_x^2)||_{L^2(\mathbb{R})} \le C_n ||v_x||_{L^\infty(\mathbb{R})}E_{n+1}^{\frac{1}{2}}E_{n+2}^{\frac{1}{2}},
\end{equation}
where for the nonlinearity term $V_3$ the Cauchy-Schwarz inequality and the Moser-type estimate \cite{Taylor_PDE} $||\partial_x^n(v_x^2)||_{L^2(\mathbb{R})}\le C_n||\partial_x v||_{L^\infty(\mathbb{R})}||\partial_x^{n+1}v||_{L^2(\mathbb{R})}$ have been used. With sufficiently small $||v_x||_{L^\infty(\mathbb{R})}$, $V_3$ can be absorbed into the two preceding terms. Later in Section  \ref{subsec: small v_x} we will prove that this quantity can be made arbitrarily small, but in this section we just assume $||v_x||_{L^\infty(\mathbb{R})}\le \sqrt{(2-1/\rho)(2n-1)}/C_n$. A constant $\rho>1/2$ can be chosen arbitrarily.

\begin{equation}\label{eq: bound of V3 using smallness of v_x}
    V_3\le \frac{1}{2}E_{n+2}+\frac{1}{2}C_n^2||v_x||_{L^\infty(\mathbb{R})}^2E_{n+1}\le \frac{1}{2}E_{n+2}+ \left(1-\frac{1}{2\rho}\right)(2n-1)E_{n+1}.
\end{equation}

By summing up all of them,

\begin{equation}\label{eq: L2 evol of E_n by higher terms}
    \frac{d}{dt}E_n \le -\frac{3}{2}E_{n+2}-\frac{2n-1}{2\rho}E_{n+1}.
\end{equation}
It shows that $E_n$ decreases with time for any $n\ge1$. To control $E_n$ with $E_{n-1}$, we adopt the interpolation inequality $E_n^2 \le E_{n-1}E_{n+1}$. 

\begin{equation}\label{eq: E_n controlled by E_n-1}
    \frac{d}{dt}\left(\frac{1}{E_n}\right)=-\frac{1}{E_n^2}\frac{d}{dt}E_n \ge \frac{2n-1}{2\rho}\frac{1}{E_{n-1}},
\end{equation}
where we neglect the negative-definite term $E_{n+2}$ of \eqref{eq: L2 evol of E_n by higher terms}.


\section{$L^\infty$ Stability of the Parabolic Steady Solution}\label{sec: stability of infinite parabola}

A linear operator $S(t)$ is well-defined due to the well-posedness of the linear perturbative equation \eqref{eq: linearized perturbed SGM} proved in Section \ref{subsec: explicit formula for linear perturbation}. According to Duhamel's principle, the solution of \eqref{eq: perturbed SGM} satisfies the following integral equation. As in the previous section, we assume the existence, uniqueness and regularities of $v$.

\begin{equation}\label{eq: Duhamel formulation}
    v(t,\cdot)=S(t)v_0-\int_0^tS(t-s)[\partial_x^2(v_x(s,\cdot) ^2)]ds.
\end{equation}

Here we try to show that the integral second term obeys a stronger bound in time, so that the overall behavior of the solution is dominated by that of the linear perturbation $w(t,x)$. For simplicity we denote $v_x(s,\cdot)^2$ as $f(s)$ and $N(s):=\partial_x^2f(s)$.

\subsection{A Priori Estimate}
In controlling the nonlinear term $N(s)$, it is important to know the behavior of the Fourier transform of it since the operator $S(t)$ involves the Fourier transform of the initial data. To do this we introduce two quantities $\Lambda(s)$ and $H_m(s)$. 

\begin{lemma}\label{lemma: a priori estimate for nonlinearity}
    For $s\ge0$ and nonnegative integer $m$, let us define $\Lambda(s)$ and $H_m(s)$ as
    \begin{equation}
        \Lambda(s):=||v_x(s,\cdot)||_{L^2(\mathbb{R})} \text{ , } H_m(s):=\left(\int_\mathbb{R}\left(1+s\xi^2\right)^m|\hat{v}(s,\xi)|^2d\xi\right)^{\frac{1}{2}}.
    \end{equation}
    Then $\hat{f}(s,k_0)$ is bounded by the following quantity.
    \begin{equation}\label{eq: boundary of f hat by Lambda and H_m}
        |\hat{f}(s,k_0)|\le\min\left(\Lambda(s)^2,\frac{\beta_ms^{-1/2}\Lambda(s)H_{m+1}(s)}{(1+sk_0^2)^{m/2}}\right).
    \end{equation}
\end{lemma}

The first part of the bound, $\Lambda(s)$, is trivial from Plancherel's theorem. For the second part, let us express the Fourier transform of $f$ by the convolution formula.

\begin{equation}\label{eq: convolution form of nonlinear fourier}
    \hat{f}(k_0)=\frac{1}{\sqrt{2\pi}}\int_\mathbb{R}\hat{v}_x(s,\xi)\hat{v}_x(s,k_0-\xi)d\xi.
\end{equation}

We split the frequency domain into two parts: $A=\{|\xi|>|k_0|/2\}$ and $B=\mathbb{R}\setminus A$. Since it always holds that either $|\xi|$ or $|k_0-\xi|$ is greater than or equal to $|k_0|/2$, the whole integral can be split into $A$ and $B$ but the integral on $B$ is bounded by that of $A$. We then apply the Cauchy-Schwarz inequality.

\begin{equation}\label{eq: second boundary of f}
\begin{split}
    |\hat{f}(k_0)|\le \sqrt{\frac{2}{\pi}}\left(\int_A |\hat{v}_x(s,\xi)|^2d\xi\right)^{1/2}\left(\int_A |\hat{v}_x(s,k_0-\xi)|^2d\xi\right)^{1/2} \le\sqrt{\frac{2}{\pi}}\Lambda(s)\left(\int_A |\hat{v}_x(s,\xi)|^2d\xi\right)^{1/2}.
\end{split}
\end{equation}

Also, we note that $(1+s\xi^2)^{m/2}\ge 2^{-m}(1+sk_0^2)^{m/2}$ in $A$. Then it follows that

\begin{equation}\label{eq: H_m boundary of f}
\begin{split}
    \left(\int_A |\hat{v}_x(s,\xi)|^2d\xi\right)^{1/2} &=\left(\int_A (1+s\xi^2)^{-m}(1+s\xi^2)^{m}|\hat{v}_x(s,\xi)|^2d\xi\right)^{1/2} \\
    &\le 2^{m}(1+sk_0^2)^{-m/2}\left(\int_A (1+s\xi^2)^{m}\xi^{2}|\hat{v}(s,\xi)|^2d\xi\right)^{1/2} \\
    &\le 2^m (1+sk_0^2)^{-m/2}s^{-1/2}H_{m+1}(s).
\end{split}
\end{equation}
Combining \eqref{eq: second boundary of f} and \eqref{eq: H_m boundary of f}, we obtain \eqref{eq: boundary of f hat by Lambda and H_m} with $\beta_m=2^{m+1/2}\pi^{-1/2}$.

\subsection{Bounds under Bootstrap Hypothesis}
As the bootstrap hypothesis, we assume $E_0(s)\le\epsilon_0^2(1+s)^{1/2}$ for $0\le s\le t$. In this section we observe how this hypothesis eventually controls the whole $v$. 

Even though we have made a hypothesis only for $E_0$, a similar estimate applies to the higher ones $E_n$. This can be done by induction, and our goal is to prove that there is a sequence $\{\epsilon_n\}$ of positive real numbers such that $E_n(t)\le \epsilon_n^2(1+t)^{1/2-n}$ for nonnegative integer $n$. Let us recall \eqref{eq: E_n controlled by E_n-1} and take $E_{n-1}(t)\le \epsilon_{n-1}^2(1+t)^{3/2-n}$, which is satisfied by $n=1$.

\begin{equation}\label{eq: E_n < (1+t)^(1/2-n)}
    \frac{1}{E_n(t)}-\frac{1}{E_n(0)} \ge \frac{1}{\rho\epsilon_{n-1}^2}\left((1+t)^{n-\frac{1}{2}}-1\right)
\end{equation}
Thus we can conclude that the bound for $E_{n-1}$ can be extended to $E_{n}\le\epsilon_n^2(1+t)^{1/2-n}$ with $\epsilon_{n}=\sqrt{\rho}\epsilon_{n-1}$, as long as \eqref{eq: bound of V3 using smallness of v_x} holds and $E_n(0)\le \rho^n\epsilon_0^2$. However, this requires $v_x$ to be small enough that $||v_x||_{L^\infty(\mathbb{R})}\le \sqrt{(2-1/\rho)(2n-1)}/C_n$. Note that if $\rho$ is close to its minimum $1/2$ then $E_n$ is bounded by a smaller quantity but the initial data should be controlled more strictly.

\subsubsection{$\Lambda(s)$ and $H_m(s)$ Controlled by the Hypothesis}

$\Lambda(s)$ is directly controlled by $E_1$.

\begin{equation}\label{eq: lambda bounded}
    \Lambda(s) = (2E_1(s))^{1/2}\le \left(2\epsilon_1^2(1+s)^{-1/2}\right)^{1/2} \le \left(2\rho\epsilon_0^2s^{-1/2}\right)^{1/2}.
\end{equation}
Next, expanding $H_m(s)$ for $0\le s\le t$, according to the binomial formula,

\begin{equation}
\begin{split}
    H_m(s)^2 &= \int_\mathbb{R}(1+s\xi^2)^m|\hat{v}(s,\xi)|^2d\xi= \int_\mathbb{R}\sum_{r=0}^{m} \binom{n}{r}(s\xi^2)^r|\hat{v}(s,\xi)|^2d\xi \\
    &= \sum_{r=0}^m \binom{m}{r}s^r \int_\mathbb{R}|\xi^r\hat{v}(s,\xi)|^2d\xi\le \sum_{r=0}^m \binom{m}{r}(1+s)^r \cdot 2E_r(s) \\
    &\le 2\left(\sum_{r=0}^{m}\rho^r\binom{m}{r}\right)\epsilon_0^2 (1+s)^{1/2}= 2(1+\rho)^m \epsilon_0^2 (1+s)^{\frac{1}{2}}.\\
\end{split}
\end{equation}
For simplicity, let us denote the coefficient $2\rho$ as $K_\Lambda$ and $2(1+\rho)^m$ as $K_{mH}$.

\begin{equation}\label{eq: H_m bounded}
    \Lambda(s)^2\le K_\Lambda\epsilon_0^2s^{-\frac{1}{2}},\  H_m(s)^2\le K_{mH} \epsilon_0^2(1+s)^{\frac{1}{2}}.
\end{equation}

\subsubsection{Time Scale of $S(t-s)N(s)$}

Now we estimate the time scale of $S(t-s)N(s)$, the integrand of the Duhamel term in \eqref{eq: Duhamel formulation}. For $0<s<1$ and $t>2$, we directly compute how the whole Duhamel term is bounded. Letting $\tau:=t-s$,

\begin{equation}
\begin{split}
    &\left|\int_0^1S(t-s)[N(s)]ds\right|\\
    &\le \int_0^1\frac{1}{\sqrt{2\pi}}\int_\mathbb{R}k_0^2\left|\hat{f}(s,k_0)\right|(1+2k_0^2\tau)^{-\frac{1}{2}}\exp\left(-\frac{k_0^4\tau}{1+2k_0^2\tau}\right)dk_0ds \\
    &\le \int_0^1\frac{1}{\sqrt{2\pi}}||\hat{v}_x(s,\cdot)||_{L^2(\mathbb{R})}^2\left(\int_{\{2k_0^2\tau<1\}}\frac{|k_0|}{\sqrt{2\tau}}dk_0+\int_{\{2k_0^2\tau\ge1\}}\frac{|k_0|\exp(-k_0^2/4)}{\sqrt{2\tau}}dk_0\right)ds \\
    &\le \frac{||v_x(0,\cdot)||_{L^2(\mathbb{R})}^2}{\sqrt{2\pi t}}\left(\frac{1}{t}+4\right)\lesssim ||v_x(0,\cdot)||_{L^2(\mathbb{R})}^2t^{-\frac{1}{2}}.\\
\end{split}
\end{equation}
The third inequality is justified by the result of Section \ref{sec: L2 estimates of the perturbation} that $||v_x(s,\cdot)||_{L^2}$ decreases in time.

In the region $s>1$, we should exploit Lemma \ref{lemma: a priori estimate for nonlinearity} and the upper bounds of $\Lambda(s)$ and $H_m(s)$. Let $I(s,\tau)$ be defined as follows, so that $|S(t-s)N(s)|\le(1/\sqrt{2\pi})I(s,t-s)$.

\begin{equation}
    I(s,\tau):=\int_\mathbb{R}k_0^2|\hat{f}(s,k_0)|(1+2k_0^2\tau)^{-1/2}\exp\left(-\frac{k_0^4\tau}{1+2k_0^2\tau}\right)dk_0.
\end{equation}

We consider the two regions of domain: the low frequency regime $\{|k_0|^2<1/s\}$ and the high frequency regime $\{|k_0|^2\ge 1/s\}$. In the low frequency regime, $|\hat{f}|$ is bounded by $\Lambda(s)^2$ and the other factors are controlled by the fact that $k_0^2s$ is small.

\begin{equation}
\begin{split}
    I_1&\le2\Lambda(s)^2\int_0^{1/\sqrt{s}}k_0^2(1+2k_0^2\tau)^{-\frac{1}{2}}\exp\left(-\frac{k_0^4\tau}{1+2k_0^2\tau}\right)dk_0 \\
    &\le2\Lambda(s)^2\int_0^{\tau/s}(\tau^{-1}\kappa)(1+2\kappa)^{-\frac{1}{2}}\left(\frac{1}{2\sqrt{\kappa\tau}}d\kappa\right) \le2^{-\frac{1}{2}}s^{-1}\tau^{-\frac{1}{2}}\Lambda(s)^2.
\end{split}
\end{equation}

Applying \eqref{eq: lambda bounded}, we obtain that $I_1\le2^{-1/2}K_\Lambda\epsilon_0^2s^{-3/2}\tau^{-1/2}$ for $s>1$. Next, in the high frequency regime, we choose the second factor of \eqref{eq: boundary of f hat by Lambda and H_m}. Also, for the proper convergence of the integral, we assume that $m>2$. It is sufficient to have that it holds for $m=3$ but Theorem \ref{thm: main result} already permits the stronger condition $m=4$.

\begin{equation}
\begin{split}
    I_2&\le \beta_m\frac{\Lambda(s)H_{m+1}(s)}{s^{1/2}}\int_{1/\sqrt{s}}^{\infty}\frac{2k_0^2(1+2k_0^2\tau)^{-1/2}}{(1+sk_0^2)^{m/2}}dk_0\le \beta_m\frac{\Lambda(s)H_{m+1}(s)}{s^{1/2}}\int_{1/\sqrt{s}}^\infty \frac{\sqrt{2}k_0^2}{(s^{m/2}k_0^m)(\tau^{1/2}k_0)}dk_0 \\
    &\le \frac{2\beta_m \sqrt{K_\Lambda K_{m+1,H}}}{m-2}\epsilon_0^2 \left(1+\frac{1}{s}\right)^{\frac{1}{2}}s^{-\frac{3}{2}}\tau^{-\frac{1}{2}}\le \frac{\beta_m\sqrt{8K_\Lambda K_{m+1,H}}}{m-2}\epsilon_0^2 s^{-\frac{3}{2}}\tau ^{-\frac{1}{2}}
\end{split}
\end{equation}
The last inequality is obtained by applying the bootstrap hypothesis. Now we have the desired regularity of $S(t-s)N(s)$.

\begin{equation}
    \left|S(t-s)N(s)\right|\le \frac{1}{\sqrt{2\pi}}\left(\sqrt{2}\rho+\frac{\beta_m \sqrt{8K_\Lambda K_{m+1,H}}}{m-2}\right)\epsilon_0^2s^{-\frac{3}{2}}(t-s)^{-\frac{1}{2}}= C_{mH}\epsilon_0^2 s^{-\frac{3}{2}}(t-s)^{-\frac{1}{2}},
\end{equation}
where in the last equality we have denoted the overall constant coefficient as $C_{mH}$. Except $s$ near zero, this results in the Duhamel term being bounded by $t^{-1/2}$.

\begin{equation}
\begin{split}
    \int_1^t(t-s)^{-\frac{1}{2}}s^{-\frac{3}{2}}ds &= \int_{\theta_0}^{\pi/2}(t-t\sin^2\theta)^{-\frac{1}{2}}(t\sin^2\theta)^{-\frac{3}{2}}(2t\sin\theta\cos\theta)d\theta \\
    &=t^{-1}\sqrt{\frac{1}{\sin^2\theta_0}-1} \lesssim t^{-\frac{1}{2}}.\\
\end{split}
\end{equation}
There has been a substitution $s=t\sin^2\theta$ and $\sin\theta_0=1/\sqrt{t}$. In conclusion, the bootstrap hypothesis $E_0(t)\le \epsilon_0^2 (1+t)^{1/2}$ leads to the $t^{-1/2}$ bound of the Duhamel term, which is strictly less than $t^{-1/2}\log t$ of the linear term $S(t)v_0$. Thus we obtain the $L^\infty$ stability stated in Theorem \ref{thm: main result}.

Now it remains for us to prove that the hypothesis keeps holding after any finite time.

\subsection{Bootstrap Closure}

The hypothesis is related to $||v(t,\cdot)||_{L^2(\mathbb{R})}$. In this section we start from the initial data satisfying $||v_0||_{L^2(\mathbb{R})}\le\epsilon_0$ and prove that $E_0(t)\le \gamma^2\epsilon_0^2(1+t)^{1/2}$ with some improving constant $\gamma<1$. 

\begin{theorem}\label{thm: bootstrap closure for small initial data}
    If we choose initial data small enough that $||v_0||_{L^2(\mathbb{R})}\le\epsilon_0$ where $\epsilon_0>0$ satisfies 
    \begin{equation}\label{eq: improving constant for bootstrap}
        \gamma:=\frac{1+K_G^{(3)}\epsilon_0}{\sqrt{2}}<1,
    \end{equation}
    and assume $E_0(s)\le \epsilon_0^2(1+s)^{1/2}$ for all $0\le s\le t$, then it follows that $E_0(t)\le \gamma^2\epsilon_0^2(1+t)^{1/2}$ at $t$. Exact formula for the constant $K_G^{(3)}$ is given in \eqref{eq: explicit K_G^(3)}.
\end{theorem}

\begin{corollary}\label{cor: E0 bootstrap closure}
    For all $t\ge0$ it holds that $E_0(t)\le\epsilon_0^2(1+t)^{1/2}$ if $||v_0||_{L^2(\mathbb{R})}$ is small enough in the sense given in Theorem \ref{thm: bootstrap closure for small initial data}, since $E_0(0)$ is strictly less than $\epsilon_0^2$.
\end{corollary}

\begin{corollary}\label{cor: high energy regularity}
    $E_n(t)\le\rho^n\epsilon_0^2(1+t)^{1/2-n}$, as long as $||v_x||_{L^\infty(\mathbb{R})}\le \sqrt{(2-1/\rho)(2n-1)}/C_n$, $E_n(0)\le\rho^n\epsilon_0^2$ and \eqref{eq: improving constant for bootstrap} hold. According to the arguments in Section \ref{prop: v_x regularity}, it is true for all $t>0$.
\end{corollary}

Let us recall the Duhamel formulation \eqref{eq: Duhamel formulation},

\begin{equation}\label{eq: boot closure - Duhamel L2}
    ||v(t,\cdot)||_{L^2(\mathbb{R})}\le ||S(t)v_0||_{L^2(\mathbb{R})}+\int_0^t||S(t-s)N(s)||_{L^2(\mathbb{R})}ds.
\end{equation}

As the first term of the right side is the solution of the linearized equation \eqref{eq: linearized perturbed SGM}, we adopt the estimate given in Lemma \ref{lemma: L2 of S(t)u}.

\begin{equation}
    ||S(t)v_0||_{L^2(\mathbb{R})}\le ||v_0||_{L^2(\mathbb{R})}(1+t)^{\frac{1}{4}}\le\epsilon_0(1+t)^{\frac{1}{4}}.
\end{equation}

Next we have to consider $||S(\tau)N(s)||_{L^2(\mathbb{R})}$. We begin with the explicit formula for $S(\tau)N(s)$, which is in the form of the inverse Fourier transform. Then we apply Plancherel's theorem to it.

\begin{equation}
\begin{split}
    ||S(\tau)N(s)||_{L^2(\mathbb{R})}^2 &= \int_\mathbb{R}(1+2k_0^2\tau)^{\frac{1}{2}}\left|k_0^2\hat{f}(s,k_0)\right|^2\exp\left(-\frac{2k_0^4\tau}{1+2k_0^2\tau}\right) dk_0 \\
    &\le \left[\sup_{k_0\in\mathbb{R}}W(k_0;s,\tau)\right] G(s). \\
\end{split}
\end{equation}
where we define $W(k_0;s,\tau):=k_0^4(1+2k_0^2\tau)^{1/2}(1+sk_0^2)^{-3}\exp(-2k_0^4\tau/(1+2k_0^2\tau))$ and $G(s):= \int(1+sk_0^2)^3|\hat{f}(s,k_0)|^2dk_0$. The factor $(1+sk_0)^3$ has been chosen so that $W$ has a proper regularity. Each term is separately considered here. When $k_0$ is low enough that $2k_0^2\tau<1$, simple calculus yields

\begin{equation}
    W(k_0;s,\tau)\le \sqrt{2}k_0^4\exp(-k_0^4\tau)\le \frac{\sqrt{2}}{e\tau}.
\end{equation}

In the remaining high frequency region, we have that $(1+2k_0^2\tau)^{1/2}\le2\sqrt{\tau}|k_0|$ and $\exp(-2k_0^4\tau/(1+2k_0^2\tau))\le \exp(-k_0^2/2)$ since $2k_0^2\tau\ge1$. 

\begin{equation}
    W(k_0;s,\tau)\le 2\sqrt{\tau}|k_0|^5(1+sk_0^2)^{-3}\exp\left(-\frac{k_0^2}{2}\right).
\end{equation}

The domain is again split into $k_0^2s\le1$ and $k_0^2s\ge1$.

\begin{equation}
    W|_{\{k_0^2s\le1\}}\le 2\tau^{\frac{1}{2}}|k_0|^5\le2\tau^{\frac{1}{2}}s^{-\frac{5}{2}},
\end{equation}

\begin{equation}
    W|_{\{k_0^2s\ge1\}}\le 2\tau^{\frac{1}{2}}s^{-3}|k_0|^{-1}\le 2\tau^{\frac{1}{2}}s^{-\frac{5}{2}}.
\end{equation}

If $s\ge1$, then $s^{-5/2}\le 2^{5/2}(1+s)^{-5/2}$ so it would not diverge after being integrated. Even though $s^{-5/2}$ is singular at $s=0$, the uniform constant upper bound of $W/\sqrt{\tau}$ justifies the bound in terms of $(1+s)$ for $0\le s\le1$ because the following inequality holds for all $k_0\in\mathbb{R}$.

\begin{equation}
    |k_0|^5\exp\left(-\frac{k_0^2}{2}\right)\le 5^{\frac{5}{2}}\exp\left(-\frac{5}{2}\right).
\end{equation}

Thus in the case of high frequency i.e. $2k_0^2\tau\ge1$, we have $W(k_0;s,\tau)\le 2\tau^{1/2}\cdot(10/e)^{5/2}(1+s)^{-5/2}$. In summary, the supremum of $W$ is estimated as

\begin{equation}\label{eq: boot closure - estimate of W}
    \sup_{k_0\in\mathbb{R}}W(k_0;s,\tau)\le \frac{\sqrt{2}}{e}\tau^{-1}+2\left(\frac{10}{e}\right)^{\frac{5}{2}}\tau^{\frac{1}{2}}(1+s)^{-\frac{5}{2}}.
\end{equation}

Next we estimate $G(s)$. The argument is very similar to what we have done to obtain $H_m(s)^2\le K_H\epsilon_0^2(1+s)^{1/2}$ under the bootstrap hypothesis. But in this case we are computing in $L^2$, so it becomes slightly different.

\begin{equation}
\begin{split}
    G(s)&= \int_\mathbb{R}(1+sk_0^2)^3|\hat{f}(s,k_0)|^2dk_0\\
    &=\sum_{r=0}^3\binom{3}{r}s^r\int_\mathbb{R} |k_0^r\hat{f}(s,k_0)|^2dk_0=\sum_{r=0}^3\binom{3}{r}s^r ||\partial_x^r(v_x^2(s,\cdot))||_{L^2(\mathbb{R})}^2.
\end{split}
\end{equation}

We distribute the derivatives to yield $\partial_x^r(v_x^2)=\sum_{j=0}^{r}\binom{r}{j}(\partial_x^{j+1}v)(\partial_x^{r-j+1}v)$. Then we apply 1D Gagliardo-Nirenberg inequality \cite{BREZIS20192839}; $||f||_{L^\infty}^2\le||f||_{L^2}||f_x||_{L^2}$. Under the bootstrap hypothesis, we are also accessible to the estimates of $E_n(s)=(1/2)||\partial_x^nv(s,\cdot)||_{L^2(\mathbb{R})}^2$.

\begin{equation}
\begin{split}
    ||\partial_x^r(v_x^2)||_{L^2(\mathbb{R})} &\le \sum_{j=0}^{r}\binom{r}{j}||\partial_x^{j+1}v||_{L^\infty(\mathbb{R})}||\partial_x^{r-j+1}v||_{L^2(\mathbb{R})} \\
    &\le \sum_{j=0}^{r}\binom{r}{j}\left(||\partial_x^{j+1}v||_{L^2(\mathbb{R})}||\partial_x^{j+2}v||_{L^2(\mathbb{R})}\right)^{\frac{1}{2}}\sqrt{2E_{r-j+1}} \\
    &\le \sum_{j=0}^{r}\binom{r}{j}2E_{j+1}^{\frac{1}{4}}E_{j+2}^{\frac{1}{4}}E_{r-j+1}^{\frac{1}{2}} \le \left[\sum_{j=0}^{r}\binom{r}{j}2\rho^{\frac{j+1}{4}}\rho^{\frac{j+2}{4}}\rho^{\frac{r-j+1}{2}}\right]\epsilon_0^2(1+s)^{-\frac{r}{2}-\frac{3}{4}} \\
    &\le 2\left[\sum_{j=0}^{r}\binom{r}{j}\rho^{\frac{2r+5}{4}}  \right]\epsilon_0^2(1+s)^{-\frac{r}{2}-\frac{3}{4}}\le 2^{r+1}\rho^{\frac{2r+5}{4}}\epsilon_0^2(1+s)^{-\frac{r}{2}-\frac{3}{4}}.
\end{split}
\end{equation}
Thus $G(s)\le K_G\epsilon_0^4(1+s)^{-3/2}$, where $K_G$ is a constant defined as

\begin{equation}\label{eq: explicit K_G}
    K_G = \sum_{r=0}^{3}\binom{3}{r}\left(2^{r+1}\rho^{\frac{2r+5}{4}}\right)^2=4\rho^{\frac{5}{2}}\sum_{r=0}^{3}\binom{3}{r}(4\rho)^r = 4\rho^{\frac{5}{2}}(1+4\rho)^3.
\end{equation}
Now we combine this with \eqref{eq: boot closure - estimate of W}.

\begin{equation}
\begin{split}
    ||S(\tau)N(s)||_{L^2(\mathbb{R})} &\le K_G^{\frac{1}{2}}\epsilon_0^2\left(\frac{\sqrt{2}}{e}\tau^{-1}(1+s)^{-\frac{3}{2}}+2\left(\frac{10}{e}\right)^{\frac{5}{2}}\tau^{\frac{1}{2}}(1+s)^{-4}\right)^{\frac{1}{2}} \\
    &\le K_G^{(1)}\epsilon_0^2\tau^{-\frac{1}{2}}(1+s)^{-\frac{3}{4}}+K_G^{(2)}\epsilon_0^2\tau^{\frac{1}{4}}(1+s)^{-2}.
\end{split}
\end{equation}
where we simply write $K_G^{(1)}=((\sqrt{2}/e)K_G)^{1/2}$ and $K_G^{(2)}=(2(10/e)^{5/2}K_G)^{1/2}$. Now we integrate $||S(t-s)N(s)||_{L^2(\mathbb{R})}$ from $s=0$ to $s=t$ to show that the right hand side of \eqref{eq: boot closure - Duhamel L2} is properly bounded. The first term is bounded by a constant;

\begin{equation}
    \int_0^t(t-s)^{-\frac{1}{2}}(1+s)^{-\frac{3}{4}}ds\le \int_0^t(t-s)^{-\frac{1}{2}}s^{-\frac{1}{2}}ds = \pi.
\end{equation}
The second term is bounded by $t^{1/4}$;

\begin{equation}
    \int_0^t(t-s)^{\frac{1}{4}}(1+s)^{-2}ds \le t^{\frac{1}{4}}\int_0^t(1+s)^{-2}ds \le t^{\frac{1}{4}}.
\end{equation}
Then $||S(\tau)N(s)||_{L^2(\mathbb{R})}$ satisfies

\begin{equation}
    \int_0^t||S(t-s)N(s)||_{L^2(\mathbb{R})}ds \le \left(\pi K_G^{(1)}+K_G^{(2)}t^{\frac{1}{4}}\right)\epsilon_0^2.
\end{equation}

Since $(a+bt^{1/4})/(1+t)^{1/4}$ is maximized at $t=(b/a)^{4/3}$ where $a,b>0$ are constants, we can get back to the form of $(1+t)^{1/4}$; that is, $\pi K_G^{(1)}+K_G^{(2)}t^{1/4}\le K_G^{(3)}(1+t)^{1/4}$. Here $K_G^{(3)}$ is a constant solely determined from the exact value of $\rho$.

\begin{equation}\label{eq: explicit K_G^(3)}
    K_G^{(3)}=\left(\left(\pi K_G^{(1)}\right)^{\frac{4}{3}}+\left(K_G^{(2)}\right)^{\frac{4}{3}}\right)^{\frac{3}{4}}
    = \left(\left(2\pi^4e^{-2}\right)^{\frac{1}{3}}+2^{\frac{2}{3}}\left(\frac{10}{e}\right)^{\frac{5}{3}}\right)^{\frac{3}{4}}K_G^{\frac{1}{2}}.
\end{equation}
Thus for all $t\ge0$,

\begin{equation}
\begin{split}
    ||v(t,\cdot)||_{L^2(\mathbb{R})}&\le \left(1+K_G^{(3)}\epsilon_0\right)\epsilon_0(1+t)^{\frac{1}{4}}, \\
    E_0(t)=\frac{1}{2}||v(t,\cdot)||_{L^2(\mathbb{R})}^2&\le\left(\frac{1+K_G^{(3)}\epsilon_0}{\sqrt{2}}\right)^2\epsilon_0^2(1+t)^{\frac{1}{2}}.
\end{split}
\end{equation}
Taking sufficiently small $\epsilon_0>0$, we finally reach the closure of the bootstrap hypothesis Theorem \ref{thm: bootstrap closure for small initial data}.


\section{Well-posedness of the Perturbative Equation}
\label{sec: well-posedness}

In the previous sections we just assumed that the perturbation $v(t,x)$ exists as a solution of \eqref{eq: perturbed SGM} and is unique with strong regularities. Now we prove the following three remaining statements. Each statement is rigorously stated in the subsection handling it.

\begin{itemize}
    \item Existence: There is a global-in-time mild solution for \eqref{eq: perturbed SGM}.
    \item Uniqueness: If $v_1$ and $v_2$ both solve \eqref{eq: perturbed SGM} with the same initial data, then they are identical.
    \item Regularity: If the initial data $v_0$ is sufficiently small, then $||v_{x}(t,\cdot)||_{L^\infty(\mathbb{R})}$ is bounded by the initial value itself.
\end{itemize}

\subsection{Proof of Existence}

The first task to be resolved is proving the existence of \eqref{eq: perturbed SGM}.

\begin{proposition}\label{prop: existence}
    For $\sigma>3/2$, let us consider a Banach space $X=C([0,T];H^\sigma(\mathbb{R}))$ equipped with the norm

    \begin{equation}
        ||f(t,\cdot)||_{X}:=\sup_{t\in[0,T]}||f(t,\cdot)||_{H^\sigma(\mathbb{R})}.
    \end{equation}
    
    If $v_0\in H^\sigma(\mathbb{R})$, then there is a time $T>0$ and a mild solution for \eqref{eq: Duhamel formulation} in $X$. Moreover, under the smallness condition for initial data mentioned in Theorem \ref{thm: main result} extended to general $\sigma$, this solution can be extended to arbitrarily large time.
\end{proposition}

In this context a mild solution means a solution of the following integral equation recalling \eqref{eq: Duhamel formulation} which is equivalent to the original \eqref{eq: perturbed SGM} if sufficient regularity is given.

\begin{equation}
    v(t,\cdot)=S(t)v_0-\int_0^tS(t-\tau)[\partial_x^2(v_x(\tau,\cdot)^2)]d\tau.
\end{equation}
For the proof of existence we adopt the contraction mapping theorem, in the closed subset $K=\{||f||_{X}\le M\}$ of $X$. Let us denote the mapping on the right hand side as $F[v]$.

\begin{equation}
    F[v]:= S(t)v_0-\int_0^tS(t-\tau)[\partial_x^2(v_x(\tau,\cdot)^2)]d\tau.
\end{equation}

We first prove that $F[K]\subseteq K$, based on the lemmata in Section \ref{sec: S(t) regularities}. For the choice of $M$ and $T$ we demand $M>2||v_0||_{H^\sigma(\mathbb{R})}$ and $T<\min\left\{1,(8C_\sigma\mu_\sigma M)^{-4}\right\}$, where $\mu_{2,\sigma,\sigma-1}$ of Lemma \ref{lemma: regularity of S(t)[D^alpha u]} has been abbreviated to $\mu_\sigma$ and $C_\sigma$ is a positive constant required for the following special case of the Kato-Ponce inequality \cite{Grfakos_Oh_Kato_Ponce_ineq}. Here $\sigma$ is required to be strictly greater than $1/2$.

\begin{equation}\label{eq: special case of Kato-Ponce}
    ||fg||_{H^\sigma(\mathbb{R})}\le C_\sigma||f||_{H^\sigma(\mathbb{R})}||g||_{H^\sigma(\mathbb{R})}.
\end{equation}
Then $||F[v]||_{H^\sigma}$ becomes less than $M$ at any time.

\begin{equation}
\begin{split}
    ||F[v]||_{H^\sigma(\mathbb{R})} &\le ||S(t)v_0||_{H^\sigma(\mathbb{R})}+\int_0^t||S(t-\tau)[\partial_x^2(v_x(\tau,\cdot)^2)]||_{H^\sigma(\mathbb{R})}d\tau \\
    &\le ||v_0||_{H^\sigma(\mathbb{R})}+\int_0^t\mu_\sigma (t-\tau)^{-\frac{3}{4}}||v_x(\tau,\cdot)^2||_{H^{\sigma-1}(\mathbb{R})}d\tau \\
    &< \frac{M}{2}+C_\sigma\mu_\sigma\sup_{t\in[0,T]}||v(t,\cdot)||_{H^{\sigma}(\mathbb{R})}^2\cdot4T^{\frac{1}{4}} \le M.
\end{split}
\end{equation}
The third inequality of the above argument is why we have demanded $\sigma>3/2$ in Proposition \ref{prop: existence}.

Secondly we check that $F$ is a contraction in $K$, which leads to the existence of $v\in X$ such that $v=F[v]$. We need to examine the nonlinear factor since $S(t)v_0$ is canceled out in $F[v_1]-F[v_2]$ where $v_1$, $v_2\in K$.

\begin{equation}
\begin{split}
    &||S(t)\partial_x^2\left[(v_{1x}^2-v_{2x}^2)(\tau,\cdot)\right]||_{H^\sigma(\mathbb{R})} \le\mu_\sigma t^{-\frac{3}{4}}||(v_{1x}^2-v_{2x}^2)(\tau,\cdot)||_{H^{\sigma-1}(\mathbb{R})} \\
    &\le C_\sigma\mu_\sigma t^{-\frac{3}{4}}||v_1+v_2||_X ||v_1-v_2||_X\le 2C_\sigma\mu_\sigma t^{-\frac{3}{4}}M||v_1-v_2||_X.\\
\end{split}
\end{equation}
Then $F$ is guaranteed to be a contraction under our condition on $T$.

\begin{equation}
    ||F[v_1]-F[v_2]||_X \le \int_0^T||S(t-\tau)\partial_x^2[(v_{1x}^2-v_{2x}^2)(\tau,\cdot)]||_Xd\tau \le 2C_\sigma\mu_\sigma M\left(4T^{\frac{1}{4}}\right)||v_1-v_2||_X. \\
\end{equation}

To show that this solution is available for all $T>0$, we adopt the results in Corollary \ref{cor: high energy regularity}. Taking $M=3||v_0||_{H^\sigma(\mathbb{R})}$ and $T=\min\{1/2,(9C_\sigma\mu_\sigma M)^{-4}\}$, we guarantee the existence of the solution in $C([0,T];H^\sigma(\mathbb{R}))$ whose initial value is $v_0$. After time $T$, we observe that $||v(T,\cdot)||_{H^\sigma(\mathbb{R})}$ is controlled by Corollary \ref{cor: high energy regularity}.
\begin{equation}
    ||v(T,\cdot)||_{H^\sigma(\mathbb{R})}^2 = \sum_{n=0}^{\sigma}\binom{\sigma}{n}||\partial_x^nv(T,\cdot)||_{L^2(\mathbb{R})}^2 \le \sum_{n=0}^\sigma \binom{\sigma}{n} 2\rho^n\epsilon_0^2(1+T)^{\frac{1}{2}-n} \le K_{\sigma H}\epsilon_0^2(1+T)^{\frac{1}{2}},
\end{equation}
where $K_{\sigma H}$ is defined in \eqref{eq: H_m bounded} as $2(1+\rho)^{\sigma}$.

If we take $v(T,x)$ as a new initial data, we can continue the solution to the time $T+T'$, and $T'$ can be given as follows.
\begin{equation}
    T'=\min\left\{\frac{1}{2},\frac{1}{\left(27C_\sigma\mu_\sigma||v_0||_{H^\sigma(\mathbb{R})}\right)^{4}(1+T)}\right\}.
\end{equation}
If the latter one is greater than $1/2$ so $T'=1/2$ is chosen, then we can repeat the same process to obtain $T$ large enough to make the latter one is chosen for the whole remaining time. If the lifetime has been lengthened to $T_{n-1}$ after $n-1$ steps, then the next longer lifetime $T_n$ is,
\begin{equation}
    T_n=T_{n-1}+\frac{c}{1+T_{n-1}}.
\end{equation}
Here the constant $c$ only depends on $\sigma$ and $v_0$. Since the sequence with this recurrence relation is divergent if the initial value is positive, the solution can be extended globally in time.

\subsection{Proof of Uniqueness}
Now let us assume $\sigma\ge5$ and $v_1$ and $v_2$ are two classical solutions of \eqref{eq: perturbed SGM} with the same initial data $v_0$. Also, let us define $\tilde{v}:=v_1-v_2$ and $V:=v_1+v_2$. We verify the uniqueness through the quantity $\int_\mathbb{R} \tilde{v}_x^2dx$.
\begin{proposition}\label{prop: uniqueness}
    $\int_\mathbb{R}\tilde{v}_x(t,x)^2dx=0$ for all $t\ge0$.
\end{proposition}

The time derivative of $\int_\mathbb{R} \tilde{v}_x^2dx$ is calculated as follows.
\begin{equation}
\begin{split}
    \frac{d}{dt}\int_\mathbb{R}\frac{1}{2}|\tilde{v}_x|^2dx &= -\int_\mathbb{R}|\tilde{v}_{xxx}|^2dx-\frac{1}{2}\int_\mathbb{R}|\tilde{v}_{xx}|^2dx + \int_\mathbb{R}\tilde{v}_{xx}\left(\tilde{v}_x V_x\right)_{xx}dx\\
    &\le -\int_\mathbb{R}|\tilde{v}_{xxx}|^2dx+\left(-\frac{1}{2}+\frac{3}{2}||V_{xx}||_{L^\infty(\mathbb{R})}\right)\int_\mathbb{R}|\tilde{v}_{xx}|^2dx+\frac{1}{2}||V_{xxxx}||_{L^\infty(\mathbb{R})}\int_\mathbb{R}|\tilde{v}_x|^2dx.\\
\end{split}
\end{equation}
We use the interpolation inequality $||\tilde{v}_{xx}||_{L^2}^2\le||\tilde{v}_x||_{L^2}||\tilde{v}_{xxx}||_{L^2}$ and $2xy\le (x^2/\lambda)+\lambda y^2$ for positive $\lambda$, $x$ and $y$. If we choose $\lambda=(1/2)|-1+3||V_{xx}||_{L^\infty(\mathbb{R})}|$, the second term is absorbed into the other terms.
\begin{equation}
\begin{split}
    \frac{d}{dt}\int_\mathbb{R}\frac{1}{2}|\tilde{v}_x|^2dx &\le -\int_\mathbb{R}|\tilde{v}_{xxx}|^2dx+\frac{1}{2}\left(\int_\mathbb{R} |\tilde{v}_{xxx}|^2dx + \lambda^2\int|\tilde{v}_x|^2dx \right)+\frac{1}{2}||V_{xxxx}||_{L^\infty(\mathbb{R})}\int_\mathbb{R}|\tilde{v}_x|^2dx\\
    &\le -\frac{1}{2}\int_\mathbb{R}|\tilde{v}_{xxx}|^2dx + \frac{1}{2}\left(\lambda^2+||V_{xxxx}||_{L^\infty(\mathbb{R})}\right)\int_\mathbb{R}|\tilde{v}_x|^2dx.
\end{split}
\end{equation}

Since the first term is negative-definite and the coefficient of the second term is uniformly bounded by quantities from initial data as $H^5(\mathbb{R})\hookrightarrow C^4(\mathbb{R})$, Gr\"onwall's lemma implies that $\tilde{v}_x$ is identically zero. Then $\tilde{v}$ must be constant in space (but might depend on time) and that constant must be zero as it is at $t=0$ and $\tilde{v}(t,\cdot)\in L^2(\mathbb{R})$.

\subsection{Uniform boundedness of $v_{x}$ in $L^\infty$}
\label{subsec: small v_x}

In extending the bootstrap hypothesis $E_0(t)\le \epsilon_0^2(1+t)^{1/2}$ to the higher order energies, we have assumed $||v_{x}(t,\cdot)||_{L^\infty(\mathbb{R})}$ is less than $\sqrt{(2-1/\rho)(2n-1)}/C_n$ at any time. Now we show that this is really an attainable property for sufficiently small initial data. This would be done by another bootstrap hypothesis; $||v_{x}(t,\cdot)||_{L^\infty(\mathbb{R})}\le \delta\le \sqrt{(2-1/\rho)(2n-1)}/C_n$ for all $t\in[0,T_\delta]$ and $1\le n\le5$.
\begin{proposition}\label{prop: v_x regularity}
    If an initial data $v_0$ of \eqref{eq: perturbed SGM} satisfies $||v_0||_{L^2(\mathbb{R})}\le\delta/\sqrt{\rho(1+\rho)}$, then $||v_{x}(t,\cdot)||_{L^\infty(\mathbb{R})}$ calculated for the solution starting from it decays in time as in \eqref{eq: v_x really decays}. Thus $||v_x(t,\cdot)||_{L^\infty(\mathbb{R})}\le\delta$ for all $t>0$.
\end{proposition}

We start from Sobolev inequality and energy estimates given by the first part of Corollary \ref{cor: high energy regularity}.
\begin{equation}
    ||v_{x}(t,\cdot)||_{L^\infty(\mathbb{R})}\le \frac{1}{\sqrt{2}}||v_{x}(t,\cdot)||_{H^1(\mathbb{R})} \le \frac{1}{\sqrt{2}}\left(||v_{x}(t,\cdot)||_{L^2(\mathbb{R})}^2+||v_{xx}(t,\cdot)||_{L^2(\mathbb{R})}^2\right)^{\frac{1}{2}} \le \left(E_1(t)+E_2(t)\right)^{\frac{1}{2}}.
\end{equation}
Since for the time $0\le t\le T_\delta$ we have assumed that $||v_{x}||_{L^\infty(\mathbb{R})}\le\delta$, we can use \eqref{eq: E_n < (1+t)^(1/2-n)}. The recurrence relation of $\epsilon_n^2$ yields $\epsilon_1^2=\rho\epsilon_0^2$ and $\epsilon_2^2=\rho^2\epsilon_0^2$.
\begin{equation}\label{eq: v_x really decays}
    ||v_{x}(t,\cdot)||_{L^\infty(\mathbb{R})} \le (\epsilon_1^2+\epsilon_2^2)^{\frac{1}{2}}(1+t)^{-\frac{1}{4}}\le\sqrt{\rho(1+\rho)}\epsilon_0(1+t)^{-\frac{1}{4}}.
\end{equation}

If $\sqrt{\rho(1+\rho)}\epsilon_0\le\delta$ then $||v_{x}||_{L^\infty(\mathbb{R})}$ becomes strictly less than $\delta$ regardless of the time $T_\delta$; hence $T_\delta$ can be extended infinitely. This also gives us the condition on $\epsilon_0$ other than what has been given in Theorem \ref{thm: bootstrap closure for small initial data}, i.e. $\gamma<1$. Combining \eqref{eq: improving constant for bootstrap}, \eqref{eq: explicit K_G} and \eqref{eq: explicit K_G^(3)}, smallness of $\epsilon_0$ needed for the stability of the parabolic steady solution is approximately
\begin{equation}
    \rho^{\frac{5}{4}}(1+4\rho)^{\frac{3}{2}}\epsilon_0<\frac{1}{40.23}.
\end{equation}
It becomes stricter as $\rho$ increases. However, for energy estimate $E_n(t)\lesssim(1+t)^{1/2-n}$, it is needed that $C_n\sqrt{\rho(1+\rho)}\epsilon_0\le\sqrt{(2-1/\rho)(2n-1)}$, which becomes small when $\rho\sim1/2$.


\appendix
\section{Classification of Steady Solutions of SGM}
In the introduction we have mentioned that there are three types of steady solutions $\tilde{u}(x)$ for \eqref{eq: SGM}. Here we describe how those solutions look like. Taking $\tilde{u}_t = 0$, every term of \eqref{eq: SGM} is at least twice differentiated. Integrating twice yields $c_1x+c_2$ on the right hand side and we restrict ourselves to the untilted case where $c_1=0$. According to the invariance under constant shift, the ODE for the steady solution becomes
\begin{equation}\label{eq: steady equation}
    \tilde{u}_{xx}+\tilde{u}+\tilde{u}_x^2=0.
\end{equation}
The following theorem presents all the possible types of steady solutions satisfying \eqref{eq: steady equation}, which will be proved throughout this appendix.

\begin{theorem}\label{thm: classification of steady solutions}
    Consider steady SGM equation \eqref{eq: steady equation} with $\tilde{u}(0)=\tilde{u}_0$ and $\tilde{u}_x(0)=0$. If $\tilde{u}_0<1/2$, then the steady solution $\tilde{u}(x)$ exists for all $x\in\mathbb{R}$ and is periodic. However, if $\tilde{u}_0>1/2$, then $\tilde{u}(x)$ blows up in negative direction in finite length so that the solution exists only for a certain bounded interval. Furthermore, the singularity at the end has logarithmic scale. Finally in the case of $\tilde{u}_0=1/2$, it has the explicit solution $u(x)=(1/2)-(1/4)x^2$ so that the solution exists in the whole real line but blows up at infinity.
\end{theorem}

\subsection{Phase Portrait Behavior}

We first analyze the phase portrait of \eqref{eq: steady equation}. Letting $\tilde{v}(x):=e^{\tilde{u}(x)}$ i.e. $\tilde{u}=\log\tilde{v}$, it is bounded below by zero as long as the solution exists.
\begin{equation}
    \tilde{v}_{xx}+\tilde{v}\log\tilde{v}=0.
\end{equation}
And this equation can be converted into a series of two first-order autonomous equations. For further discussions we would treat the independent variable $x$ as the time $t$ in dynamical systems theory.
\begin{equation}\label{steady ODE in 2D 1st order}
    \left\{ \begin{array}{ccc}
    p'(t) &=& q(t)  \\
    q'(t) &=&  -p(t)\log p(t) \\
    (p(0),q(0)) &=& (p_0,q_0)
     \end{array}
    \right.
\end{equation}

Note that $p\ge0$ with the convention $p\log p|_{p=0}=0$. Then $(p,q)=(1,0)$ is the only fixed point of the system. Defining $F(p,q):=(q,-p\log p)$, then $F\in C^{\infty}((0,\infty)\times\mathbb{R})$ and eigenvalues of $\nabla F(1,0)$ are $\pm i$. Since this system has the Hamiltonian $H(p,q)=\frac{1}{4}p^2(1-2\log p)-\frac{1}{2}q^2$ conserved throughout the trajectory, level sets near $(1,0)$ become closed curves. As $H(1,0)=1/4$ and $H(0,q)=-\frac{1}{2}q^2<1/4 = H(1,0)$, closed orbit defined by $H(p,q)=E_0$ with $E_0\in(0,1/4)$ never touches the $q$-axis. Furthermore, note that $H(p^*,0)=0$ is satisfied by $p^*=0$ and $p^*=\sqrt{e}$, which corresponds to $\tilde{u}_0=1/2$.

\begin{figure}[htbp]
\centering

\begin{subfigure}{0.49\textwidth}
\includegraphics[width=\linewidth]{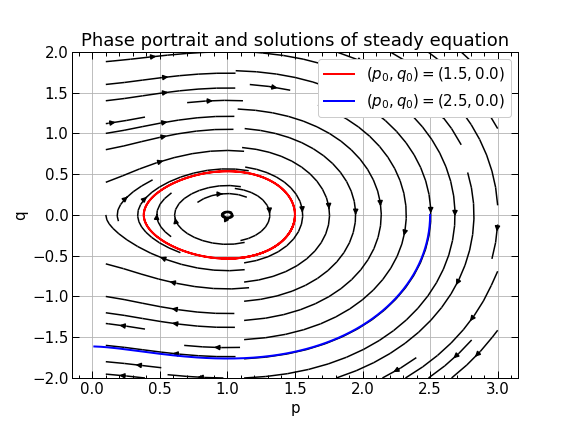}
\caption{Phase portrait of \eqref{steady ODE in 2D 1st order}. Red closed orbit is the solution starting from $(1.5,0.0)$ and blue orbit starts from $(2.5,0.0)$. One can see that the red orbit is periodic while the blue is not and ends up touching the $q$-axis, which represents the negative blow-up of the steady solution $\tilde{u}(x)$.}
\label{fig: phase portrait and solutions of steady equation}
\end{subfigure}
\hfill
\begin{subfigure}{0.49\textwidth}
\includegraphics[width=\linewidth]{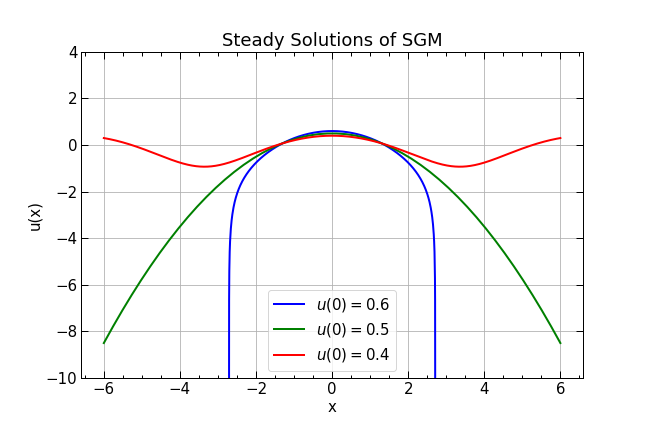}
\caption{This figure demonstrates the different types of solutions classified by Theorem \ref{thm: classification of steady solutions}. Blue line represents the solution from $\tilde{u}(0)=0.6$ and red line for $\tilde{u}(0)=0.4$. $\tilde{u}(0)>0.5$ leads to finite length blow-up and $\tilde{u}(0)<0.5$ leads to periodic orbit. Green line represents the critical case $\tilde{u}(0)=0.5$.}
\label{fig: Steady solution comparison: periodic and diverging}
\end{subfigure}

\caption{Visualizations of steady solutions, in phase portrait and in real space.}
\end{figure}

Meanwhile, if the constant value of Hamiltonian $E_0$ is negative, $q(t)<0$ for the whole lifetime of the solution because $q$ consistently decreases and becomes negative while passing through the region $\Omega_1:=\{p\ge1,\ q\le0\}$ but it cannot touch zero while in the region $\Omega_2:=\{0<p<1,\ q\le0\}$. The solution from $\Omega_1$ is guaranteed to get into $\Omega_2$ since $p'(t)=q(t)<0$ in this region and its absolute value keeps increasing as time goes on. Recalling that $H(p,q)=\frac{1}{4}p^2(1-2\log p)-\frac{1}{2}q^2=E_0<0$ and exploiting the fact that $p^2(1-2\log p)>0$ for $p\in(0,\sqrt{e})$, we can see that $q<-\sqrt{-2E_0}$. Combining with $p'(t)=q(t)$, we obtain the result that the solution blows up after some time less than $\sqrt{e/(-2E_0)}$. Also, the same process for $q'=-p\log p$ and $p>\sqrt{e}$ concludes that it takes finite time from the initial point to the vertical line $p=\sqrt{e}$. Gathering up all of them, when the solution starts from initial data $(p_0,0)$ where $p_0>\sqrt{e}$, there exists a finite positive real number $T^*$ such that $p(T^*)=0$. It means that there exists a finite positive real number $L$ such that $\lim_{x\uparrow L}\tilde{u}(x)=-\infty$.

As the final case we investigate what happens when $(p_0,q_0)=(\sqrt{e},0)$ i.e. $\tilde{u}(0)=1/2$. In this case, $\tilde{u}(x)=(1/2)-(1/4)x^2$ works as the explicit solution. Due to the uniqueness of the system \eqref{steady ODE in 2D 1st order}, this must be the only solution satisfying \eqref{eq: steady equation}. Note that the original time-dependent equation \eqref{eq: SGM} is invariant under translation, so that $u(t,x)=a-(1/4)(x-b)^2$ is a steady solution for \eqref{eq: SGM}.

\subsection{Singularity Scale of Steady Solutions}
As we know the condition under which the steady solution forms singularities in finite length, now we investigate how those singularities look like. Note that the singularity near $x\sim L$ corresponds to $p\sim 0$ or $t\sim T^*$. In this region $dq/dp\sim0$ along the solution.
\begin{equation}
    \frac{dq}{dp} = \frac{q'}{p'} = -\frac{p\log p}{q}\sim0
\end{equation}

Since the orbit is almost parallel to the $p$-axis and the velocity is $(p',q')\sim(q,0)$, $q\sim-\bar{q}$ for some constant $\bar{q}>0$ and $p(t)\sim\bar{q}(T^*-t)$. Also from the definition $p=\tilde{v}=e^{\tilde{u}}$ we obtain $\tilde{u}(x)=\log p(x)\sim\log(L-x)$. This is the same for the opposite side due to the symmetry of the equation under $x\rightarrow-x$. Therefore the solution exists in the bounded interval $(-L,L)$ and it blows up in logarithmic scale.

\section*{Acknowledgements}
The author is deeply grateful to In-Jee Jeong for a number of valuable discussions and comments throughout this paper. This work was supported by the grant RS-2024-00406821.

\nocite{*}
\bibliography{Stable_Inverted_Parabola_bibtex}
\bibliographystyle{amsalpha}

\end{document}